\documentclass[final,onefignum,onetabnum]{siamart171218}
\usepackage[T1]{fontenc}

\usepackage{listings}
\usepackage{lipsum}
\usepackage{amsfonts}
\usepackage{graphicx}
\usepackage{epstopdf}
\usepackage{algorithmic}
\ifpdf
  \DeclareGraphicsExtensions{.eps,.pdf,.png,.jpg}
\else
  \DeclareGraphicsExtensions{.eps}
\fi

\newsiamremark{remark}{Remark}
\newsiamremark{hypothesis}{Hypothesis}
\crefname{hypothesis}{Hypothesis}{Hypotheses}
\newsiamthm{claim}{Claim}

\headers{Parallel Simulations of Biofilm Flow}{N. McClanahan, N. Stegmeier, R. Sundermann J.  Doom and J-H. Kimn}

\title{Parallel Simulations of Biofilm Flow using the Modified Cahn-Hilliard Equation\thanks{
\funding{This material is based upon work supported by the National Science Foundation under Grant No. 1559978.}}}

\author{Nathan McClanahan\thanks{Department of Mathematics and Statistics, South Dakota State University, Brookings, SD (\email{Nathan.McClanahan@sdstate.edu}).}
\and Nicholas Stegmeier\thanks{Department of Mathematics and Statistics, South Dakota State University, Brookings, SD (\email{Nicholas.Stegmeier@sdstate.edu}).}
\and Rylee Sundermann\thanks{Department of Mathematics and Statistics, South Dakota State University, Brookings, SD (\email{Rylee.Sundermann@sdstate.edu}).}
\and Jeffrey Doom\thanks{Department of Mechenical Engineering, South Dakota State University, Brookings, SD (\email{Jeffrey.Doom@sdstate.edu}).}
\and Jung-Han Kimn\thanks{Department of Mathematics and Statistics, South Dakota State University, Brookings, SD (\email{Jung-Han.Kimn@sdstate.edu}).}}

\usepackage{amsopn}

\ifpdf
\hypersetup{
  pdftitle={Parallel Simulations of Biofilm Flow using the Modified Cahn-Hilliard Equation},
  pdfauthor={Nathan McClanahan, Nicholas Stegmeier, Rylee Sundermann, Jeffrey Doom, and Jung-Han Kimn}
}
\fi

\begin{document}

\maketitle

\begin{abstract}
We present a 2D parallel implementation of the modified Cahn-Hilliard equation for the simulation of  a  biofilm in an aqueous enviroment. Biofilms are attached microbial communities made of many different components and can have both positive and negative effects. They can be used for bioremediation but also are the cause of the majority of chronic infections. It is for these reasons that we study them. Due to their composition being mostly water we choose to model them as an incompressible fluid. We used a visco-elastic phase separation model based on the modified Cahn-Hilliard equation and the Flory-Huggins energy density. We present results of a simulation showing detachment of a biofilm protrusion from a base layer of biofilm due to the flow over the biofilm. The parallelization was accomplished using PETSc (Portable, Extensible Toolkit for Scientific Computing), specifically the DMDA abstraction layer within PETSc. DMDA provides a useful interface for the solution of linear systems arising from structured grid discretizations. We evaluate the parallel performance of the implementation with a strong scaling test and calculate the speedup for various mesh sizes. 

\end{abstract}



\section{Introduction}
\label{sec:Intro}

Biofilms are found throughout the world.  They are common in both nature and in man made enviroments.  Biofilms are a collection of micro-organisms that adhere to a surface through a self-produced extracellular polymeric substance (EPS).  This is a sticky, slimy stubstance that holds the biofilm in place. Biofilms cause massive losses around the world.  They cost the U.S. alone billions of dollars every year.  Biofilms contaminate food, water, and industrial equipment.  They are also estimated to be the leading cause of chronic infections in the world\cite{costerton2003medical}.  Biofilms are not all bad though.  They play a role in bioremediation, they are a safer, less toxic way to mine certain hard to obtain minerals, and they can be used as sealants through biomineralization\cite{tabak2005developments, phillips2013engineered, spath1998sorption}. \par
Due to widespread nature of biofilms, the massive monetary losses caused by them, and the potential benefits from using them the last three decades have seen an increased desire to study them\cite{picioreanu1998mathematical, picioreanu1999multidimensional, picioreanu2004particle, picioreanu2004advances}.  There are models that are discrete where the biofilm is modeled using a cellular automaton approach\cite{picioreanu1996modelling} and there are models that are continuum based\cite{alpkvista2007multidimensional}.  There are also models that are a combination of the two, where the cellular automaton approach is used for the biofilm and the continuum model is used for the nutrient substrate\cite{picioreanu1998new}.  \par
There are many more aspects of a model besides just discrete or continuum based.  There are single fluid and multifluid models\cite{cogan2005channel, Zhang_1D}.  There are models with multiple compenents to the fluid\cite{Zhang_1D}.  Varying the number of species of organisms in the biofilm can produce different results\cite{alpkvista2007multidimensional}.  Biofilms have been treated as both elastic and visco-elastic fluids depending on the time scale used\cite{Zhang_1D}.  Experimental results support the use of a visco-elastic fluid model\cite{klapper2002viscoelastic}. \par

High fidelity simulation of biofilms is numerically challenging and computationally expensive. In addition to traditional concerns from computational fluid dynamics such as ensuring a divergence free velocity field and coupling the pressure and velocity, the Cahn-Hilliard equation introduces a nonlinear, fourth order equation coupled to the momentum equation. This necessitates a larger stencil and more interprocessor communication. The discretized equations result in several sparse matrix systems of the form $Ax=b$. Direct solution of these systems is often intractable due to their large size, so parallel iterative methods are pursued instead. For the parallelization, we use the Portable, Extensible Toolkit for Scientific Computation (PETSc) \cite{petsc-web-page, petsc-user-ref} to decompose the computational domain and iteratively solve the resulting linear systems. In particular, the DMDA abstraction layer in PETSc is used as an efficient interface for mapping the computational domain to the matrix.  There has been previous work done using parallel computational methods and variations of the Cahn-Hilliard equation.  One of these include work done by Zheng, Yang, Cai, and Keyes \cite{zheng2015parallel}.  Their model differs in that they use the Cahn-Hilliard-Cook equation instead which has an additional term to model noise, in this case thermal fluctuations, in the system.

This paper is organized as follows. We review the mathematical models developed by Zhang et al \cite{Zhang_1D, Zhang_2D} in \cref{sec:MathModel} and \cref{sec:NonDim}. The parallel implementation and numerical results are described in \cref{sec:NumericMethod} and \cref{sec:NumericResults}, and the conclusions follow in
\cref{sec:Conc}.

\section{Mathematical Model}
\label{sec:MathModel}
We follow the one-fluid two-compenent model developed in \cite{Zhang_1D, Zhang_2D}.  We define the following variables.  Let the average velocity be $\mathbf{v}$, the pressure be $p$, the volume fraction of the polymer be $\phi_n$, the volume fraction of the solvent be $\phi_s$, and $c$ be the nutrient concentration.
\subsection{Momentum Transport and Continuity}
We assume incompressible flow which results in the following continuity equations for momentum transport 
\begin{align}
\nabla\cdot\mathbf{v}&=0,\label{eq:divfree1d}\\
\rho\frac{d\mathbf{v}}{dt}& = \nabla\cdot\left(\phi_n\tau_n+\phi_s\tau_s\right)-\left[\nabla p +\gamma_1kT\nabla \cdot\left(\nabla\phi_n\nabla\phi_n\right)\right], \label{eq:momentumtrans}
\end{align}
where the density is a volume fraction averaged density $\rho = \phi_n\rho_n+\phi_s\rho_s$ with $\rho_n$ being the density of the polymer network, $\rho_s$ is the density of solvent, $\tau_n$ is the extra stress tensor for the network, $\tau_s$ is the extra stress tensor for the solvent, $k$ is the Boltzman constant, $T$ is the temperature, and $\gamma_1$ measures the strength of the conformation entropy.  We use a volume fraction averaged velocity field $\mathbf{v}=\phi_n\mathbf{v}_n+\phi_s\mathbf{v}_s$.  We use a phase separation energy functional instead of the extended Flory-Huggin's mixing free energy density for ease of use
\begin{equation}\label{eq:modenergy}
f(\phi_n) = kT\left[\frac{\gamma_1}{2}||\nabla\phi_n||^2 + \gamma_2kT\left(\phi_n^2\left(1-\phi_n\right)^2\right)\right].
\end{equation}
where $\gamma_2$ is the strength of the bulk free energy from the Flory-Huggin's free energy density.  The incompressible condition also implies that 
\begin{equation}
\phi_n+\phi_s=1\notag
\end{equation}

\subsection{The Cahn-Hilliard Equation as the Network Transport Equation}\label{sec:CHeq}
We used the singular or modified Cahn-Hilliard equation as the polymer network transport equation
\begin{equation}
\frac{\partial \phi_n}{\partial t} + \nabla \cdot \left(\phi_n\mathbf{v}\right)=\nabla\cdot\left(\lambda\phi_n\nabla\frac{\delta f}{\delta \phi}\right)+g_n.\label{eq:modCH}
\end{equation} 
where $\lambda$ is the mobility parameter and the network production rate is given by
\begin{equation}
g_n = \epsilon\mu \phi_n \frac{c}{K_c+c'}
\end{equation}
where $\epsilon$ is a scaling constant, $\mu$ is maximum production rate, and $K_c$ is the half-sturation constant. The Cahn-Hilliard equation has a polymer network volume fraction dependent mobility.  We chose to use this equation as opposed to the standard Cahn-Hilliard equation as it has been shown that the modified Cahn-Hilliard equation is more appropriate to use when modeling the transport of the polymer network especially with production of polymer included\cite{Zhang_1D}.
\subsection{Nutrient Transport Equations}
We use the nutrient transport equation
\begin{equation}
\frac{\partial}{\partial t}\left(\phi_sc\right) + \nabla\cdot\left(c\mathbf{v}\phi_s-D_s\phi_s\nabla c\right) = -g_c
\end{equation}
where $c$ is the nutrient concentration, $D_s$ is the diffusion coefficient for the nutrient substrate, and the nutrient consumption rate is given by 
\begin{equation}
g_c = \phi_nAc
\end{equation}
with $A$ being a constant.

\subsection{Constitutive Equations}The constituive equations are as follows
\begin{equation}\label{eq:tau}
\tau_n=2\eta_n\mathbf{D}, \hspace{.25 in} \tau_s =2\eta_s \mathbf{D},
\end{equation}
where $\eta_n$ and $\eta_s$ are the viscosities of the polymer network and solvent.  The rate of strain tensor, $\mathbf{D}$, is given by 
\begin{equation}\label{eq:D}
\mathbf{D}=\frac{1}{2}\left[\nabla\mathbf{v}+\nabla \mathbf{v}^T\right].
\end{equation}
The velocities $\mathbf{v}_n$ and $\mathbf{v}_s$ have two parts to them.  The first part is the convection due to the average velocity $\mathbf{v}$.  For $\mathbf{v}_n$ the second part is an excessive flux due to the mixing of the polymer network and solvent.  The polymer network excessive flux is defined as proportional to the gradient of the variation of the free energy.
\begin{equation}
\mathbf{v}^e_n = -\lambda \phi_n \nabla \frac{\delta f}{\delta \phi_n}
\end{equation}
The excessive flux for the solvent is due to the spatial inhomogeneity of the mixture and is defined as \cite{Zhang_1D}.
\begin{equation}
\mathbf{v}^e_s = \frac{\lambda \phi_n}{\phi_s}\nabla \frac{\delta f}{\delta \phi_n}
\end{equation}
 We define $\mathbf{v}_n$ and $\mathbf{v}_s$ as follows 
\begin{equation}
\mathbf{v}_n=\mathbf{v}+\mathbf{v}^e_n \hspace{.4 in}
\mathbf{v}_s=\mathbf{v}+\mathbf{v}_s^e.
\end{equation}
For the cases where shear flow is present at the top boundary we use periodic boundary conditions in the $x$ direction and the following boundary conditions in the $y$ direction.
\begin{align}\label{eq:BCanalytic}
&\nabla \left[c\mathbf{v}\phi_s -D_s\phi_s\nabla c\right]\cdot\mathbf{n}\big{|}_{y=0,1}=0, \hspace{.2 in} \nabla \phi_n\cdot\mathbf{n}\big{|}_{y=0,1}=0,\notag\\
&\nabla \left[\mathbf{v}\phi_n -\Lambda\phi_n\frac{\delta f}{\delta \phi_n}\right]\cdot\mathbf{n}\big{|}_{y=0,1}=0, \hspace{.25 in} \mathbf{v}\big{|}_{y=0}=\mathbf{0}, \hspace{.2 in} \mathbf{v}\big{|}_{y=1}=\mathbf{v}_0
\end{align}

\section{Nondimensionalization}
\label{sec:NonDim}
The system of equations will be nondimensionalized using a characteristic time-scale, $t_0$, and length-scale, $h$.  The values are specified in a table below in the results section.  The nondimensionalized variables are as follows:
\begin{equation} \label{eq:nondimvars3}
\widetilde{t} = \frac{t}{t_0}, \hspace{.25 in} \widetilde{\mathbf{x}} = \frac{\mathbf{x}}{h}, \hspace{0.25 in} \widetilde{\mathbf{v}} = \frac{\mathbf{v}t_0}{h}, \hspace{0.25 in} \widetilde{p}=\frac{pt_0^2}{\rho_0h^2}, \hspace{0.25 in} \widetilde{\tau}_n=\frac{\tau_nt_0^2}{\rho_0h^2}, \hspace{0.25 in} \widetilde{c}=\frac{c}{c_0}
\end{equation}
with $c_0$ being a characteristic substrate concentration.  Using these variables the following nondimensional parameters are found\cite{Zhang_1D}
\begin{align} \label{eq:nondimparams3}
\hspace{.2 in}\Lambda =& \frac{\lambda p_0}{t_0}, \hspace{.25 in} \Gamma_1=\frac{\gamma_1kTt_0^2}{\rho_0h^4},\hspace{.25 in} \Gamma_2 = \frac{\gamma_2kTt_0^2}{\rho_0h^2},\notag \\
Re_s=&\frac{\rho_0h^2}{\eta_st_0},\hspace{.25 in} Re_n=\frac{\rho_0h^2}{\eta_nt_0},\hspace{.25 in}  
\widetilde{D}_s=\frac{D_st_0}{h^2},\\
\widetilde{\rho}=&\phi_s\frac{\rho_s}{\rho_0}+\phi_n\frac{\rho_n}{\rho_0},\hspace{.15 in} \widetilde{A}=At_0,\hspace{.15 in} \widetilde{\mu}=\mu t_0\notag 
\end{align}
where $\rho_0$ is the averaged density; $Re_{s,n,p}$ are the Reynolds numbers for the solvent and the polymer network.  The $\Lambda, \, \Gamma_{1,2},\, \widetilde{D}_s,$ and $\widetilde{\mu}$ are the dimensionless versions of the same named dimensional parameters.  The Deborah number, $\Lambda_1$, is a number that is used to characterize the fluidity of materials\cite{reiner1964deborah}.  In this case it is used in the dimensionless version of the differential equation to solve for $\tau_n$ from equation \eqref{eq:tau}.\par
For simplicity we drop the $~\widetilde{}$ and the dimensionless equations are then
\begin{align}\label{eq:dimless}
&\nabla\cdot \mathbf{v}=0,\notag \\
&\rho\frac{d\mathbf{v}}{dt}=\nabla \cdot\left(\phi_n\tau_n+\phi_s\tau_s\right) - \left[\nabla p +\Gamma_1\left(\nabla\phi_n\nabla\phi_n\right)\right],\notag\\
&\frac{\partial}{\partial t}\left(\phi_s c\right) + \nabla \cdot\left(c\mathbf{v}\phi_s-D_s\phi_s\nabla c\right) = -g_c,\\
&\frac{\partial \phi_n}{\partial t} +  \nabla  \cdot \left(\phi_n\mathbf{v}\right)=\nabla\cdot\left(\Lambda\phi_n\nabla\frac{\delta f}{\delta \phi}\right)+g_n.\notag\\
&f(\phi_n) = \frac{\Gamma_1}{2}||\nabla\phi_n||^2 + \Gamma_2\left(\phi_n^2\left(1-\phi_n\right)^2\right)\notag
\end{align}
\label{sec:nondim}

\section{Numerical Methods}
\label{sec:NumericMethod}
Looking to equation \ref{eq:dimless} as a starting point we first discuss the continuity and momentum transport equations.  We solve the momentum transport and continuity equations using a velocity corrected projection method.  The projection method first proposed by Alexandre Chorin in 1967 works equally well for both 3D and 2D\cite{chorin1967}.  This method decouples the equations so they are easier to solve.  It only requires solving two decoupled equations for pressure and velocity which makes it efficient for numerical simulations\cite{guermond_projectionmethod}.\par
The projection method proposed by Chorin uses the incompressible Navier-Stokes equation as a starting point.  Our momentum equation needs to be re-formulated in order to fit this method.  Let 
\begin{equation}
\mathbf{R}= -\nabla \cdot \left( \Gamma_1 \nabla \phi_n \nabla \phi_n \right) + \nabla \cdot \left( \phi_n \tau_n + \phi_s \tau_s- \frac{2}{Re_a}\mathbf{D} \right)
\end{equation}
where $\mathbf{D}$ is as defined in equation \ref{eq:D} and $Re_a$ is the volume fraction averaged Reynolds number.  This then gives rise to the momentum transport equation as 
\begin{equation}\label{eq:momtoNS}
\rho \left( \frac{\partial \mathbf{v}}{\partial t}+ \mathbf{v}\cdot \nabla \mathbf{v}\right) =-\nabla p +\frac{1}{Re_a} \nabla ^2 \mathbf{v}+\mathbf{R} .
\end{equation}
Using the projection method we can solve for $\mathbf{v}$ using the following process.  If we ignore the pressure in equation \ref{eq:momtoNS} we can solve the following boundary value problem for $\mathbf{u^{n+1}}$. 
\begin{align}\label{eq:nopressmomentum}
\frac{\mathbf{u^{n+1}}-\mathbf{v}^n}{\Delta t}=-\mathbf{v}^n\cdot \nabla \mathbf{v}^n +\frac{1}{\rho Re_a}\nabla ^2 \mathbf{v}\\
\mathbf{u}^{n+1}\Big{|}_{y=0}=0,\hspace{.2 in} \mathbf{u}^{n+1} \Big{|}_{y=1} = \mathbf{u}_0\notag
\end{align}
This will result in a velocity, $\mathbf{u^{n+1}}$, that is not divergence free.  Next we solve for the pressure at the next time step so that we can use that to correct the velocity.  We solve the following Poisson equation with Neumann boundary condition for the pressure.
\begin{align}\label{eq:pressuresolve}
&-\nabla \cdot \frac{1}{\rho^{n+1}}\nabla p^{n+1} = \nabla \cdot \mathbf{u}^{n+1}\\
&\frac{\partial p^{n+1}}{\partial n}\bigg{|}_{y=0,1}=0\notag
\end{align}
With this updated pressure we update the velocity to enforce the divergence free condition using equation \ref{eq:velocityupdate}. 
\begin{equation}\label{eq:velocityupdate}
\mathbf{v}^{n+1} = \mathbf{u}^{n+1}+\frac{1}{\rho^{n+1}}\nabla p^{n+1}
\end{equation}
We used central differences for the spatial discretization to ensure second order accuracy.  In order to get second order accuracy in time we use Crank-Nicholson and extrapolation for the non-linear terms in $\mathbf{R}$ and $f$.  We used a structured grid with uniform mesh size in both space and time.  The time step size is denoted using $\Delta t$ and the spatial step sizes are denoted as $\Delta x$ and $\Delta y$.  The computational domain $\Omega = [0,1]\times [0,1]$ is divided using nodes located at $(x_i,y_j) = (i\Delta x, j \Delta y)$ with $i = 0, 1, ... N_x$ and $j = 0,1,...,N_y$.  At the node points $(n\Delta t, i\Delta x, j\Delta y)$ we denote the solutions using a superscript for time and subscripts for space.  For the network volume fraction we use $\phi^n_{n,i,j}$ and for the nutrient concentration we use $c^n_{i,j}$.  For the cases covered in the results section \ref{sec:NumericResults} we use the boundary condition $\mathbf{v\cdot n}\big{|}_{y=0,1} =0$.  This results in the boundary conditions for $\phi_n$ and $c$ in equation \ref{eq:BCanalytic} becoming 
\begin{equation}
\nabla c \cdot \mathbf{n}\big{|}_{y=0,1}=0, \hspace{.2 in} \nabla \phi_n\cdot \mathbf{n}\big{|}_{y=0,1} =0,\hspace{.2 in} \nabla \frac{\delta f}{\delta \phi_n}\cdot \mathbf{n}\big{|}_{y=0,1}=0.
\end{equation}
This would represent a zero flux of these quantities through the corresponding surface.  These boundary conditions result in the discrete boundary conditions given by the following equations.\\
\begin{align}\label{eq:discreteBC}
&\phi^n_{n,i,1}=\phi^n_{n,i,-1}, \hspace{.2 in}  \phi^n_{n,i,2} = \phi^n_{n,i,-2},\notag \\
&\phi^n_{n,i,N_{y}+1}=\phi^n_{n,i,N_{y}-1}, \hspace{.2 in}  \phi^n_{n,i,N_{y}+2} = \phi^n_{n,i,N_y-2}, \\
&c^n_{i,1}=c^n_{i,-1}, \hspace{.2 in}  c^n_{n,i,N_{y}+1} = c^n_{n,i,N_{y}-1}\notag
\end{align}

\subsection{Parallel Implementation} 

Discretization of system \ref{eq:dimless} results in a sequence of linear algebraic systems. For high fidelity simulations, numerical solution of these systems is intractable without the use of parallel computing. We use PETSc and MPI (the Message Passing Interface) to parallelize the data structures and iteratively solve the linear systems. This is accomplished through the Data Management for Distributed Arrays (DMDA) object in PETSc. DMDA provides functionality to decompose the computational domain, map domain indices to matrix indices for efficient matrix construction, and can simplify the implementation of boundary conditions. In conjunction with PETSc, MPI is used to communicate ghost points between processors. In short, we use PETSc to handle the domain-to-matrix mapping and matrix computation, while MPI is used to communicate multi-dimensional array data. This approach combines the efficient, scalable numerical methods in PETSc with the flexibility and control of parallel communication with MPI. 

To provide a concrete example of how we use PETSc, DMDA, and MPI, consider a computational domain of \texttt{Nx=4} by \texttt{Ny=4} nodes.
\begin{figure}[h!]
	\centering
	\includegraphics[width=0.5\linewidth,trim={0 35 0 14.75}]{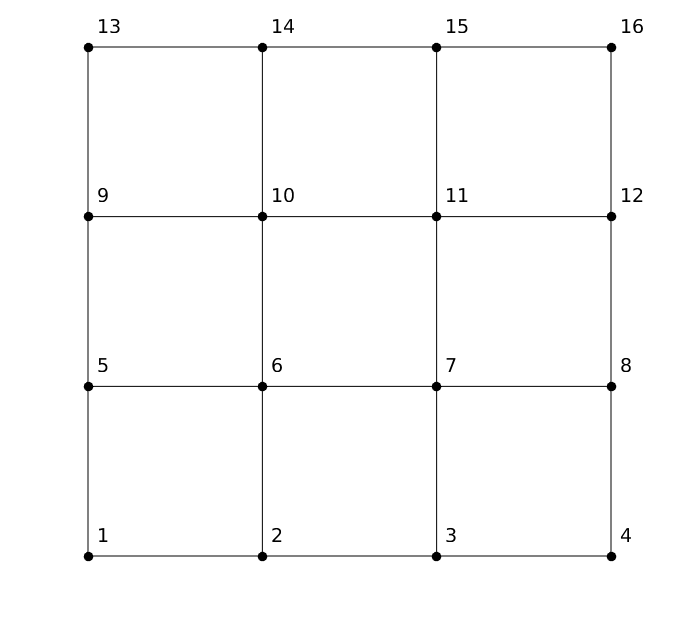}
	\caption{An example computational domain when \texttt{Nx=4}, \texttt{Ny=4}.}
	\label{fig:axb}
\end{figure}
PETSc's DMDA object must be provided the number of nodal points in each direction, the stencil type and width, and the boundary conditions in each direction. Given this information, PETSc determines an efficient way of partitioning the mesh among processors. The DMDA object is initialized using \texttt{DMDACreate2D()}.

\begin{lstlisting}[label={lst:dmda_create},caption={Initializing a 2D DMDA object in PETSC.}]
// Inputs:			Description:
// PETSC_COMM_WORLD		PETSc communicator
// bx,by			boundary types			
// stype			stencil type
// Nx,Ny		    	number of nodes	
// dof				degrees of freedom
// sw				stencil width
// da				DMDA object		
DMDACreate2d(PETSC_COMM_WORLD,bx,by,stype,Nx,Ny,
PETSC_DECIDE,PETSC_DECIDE,dof,sw,NULL,NULL,&da);
\end{lstlisting}

The details of the resulting decomposition will depend on the number of nodes and number of processors in use. For a two processor case, DMDA might decompose the domain as follows:

\begin{figure}[ht] 
	\label{ fig7} 
	\begin{minipage}[tb]{0.5\linewidth}
		\centering
		\includegraphics[width=.7\linewidth]{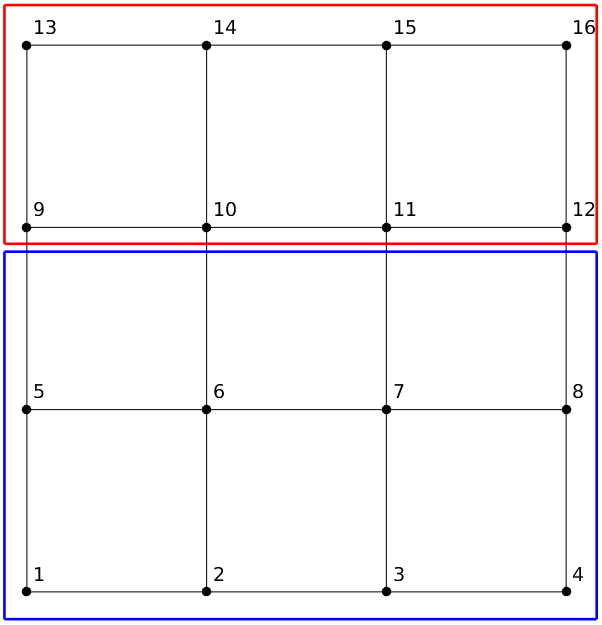} 
		\vspace{4ex}
	\end{minipage}
	\begin{minipage}[tb]{0.5\linewidth}
		\centering
		\includegraphics[width=.8\linewidth]{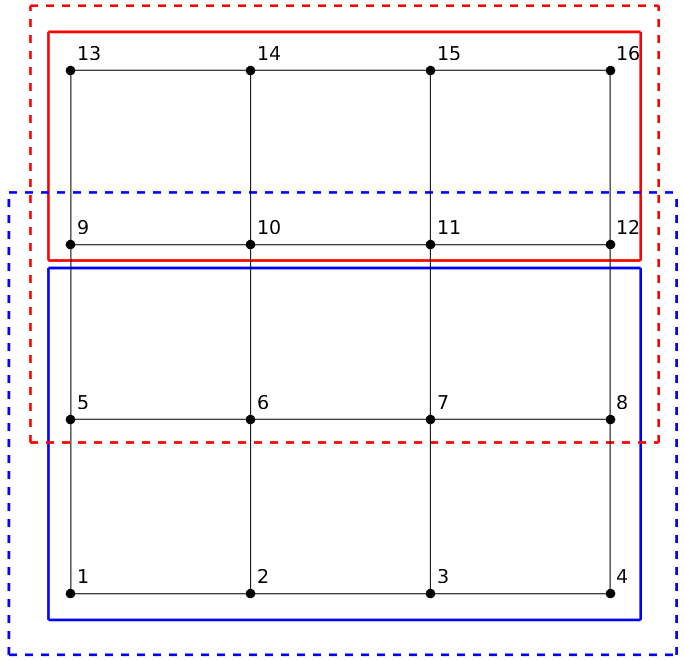} 
		\vspace{4ex}
	\end{minipage}
	\caption{DMDA partitioning of mesh for two processors, with and without ghost points.}
\end{figure}

One of the more useful features of DMDA is the ability to associate matrices and vectors with the domain decomposition. A matrix \texttt{A}, solution vector \texttt{x} and right-hand side vector \texttt{b} can be associated with the \texttt{da} object created in Listing~\ref{lst:dmda_create}.

\begin{lstlisting}[label={lst:dmda_matvec},caption={Associating matrices and vectors with DMDA.}]
// Create PETSc matrix A and vectors x,b
DMCreateMatrix(da,&A);
DMCreateGlobalVector(da,&x);
DMCreateGlobalVector(da,&b);
\end{lstlisting}

Normally, the user must manually decompose and parallelize the computational domain and then paralellize matrices and vectors in a manner that is conformal with the domain parallelization. In Listing~\ref{lst:dmda_matvec}, this process is handled automatically. Furthermore, when constructing the matrix it is usually necessary to carefully map the domain indices to matrix rows and columns. PETSc offers the \texttt{MatStencil} data structure as way to automatically compute these index transformations. For example, at node a \texttt{(i,j)} in the interior of the domain, a matrix row corresponding to the discrete 2D Laplacian could be entered as follows:
\begin{lstlisting}[label={lst:dmda_entry},caption={Constructing the matrix with DMDA.}]
// Value and index data structures
PetscReal	   	vals[5];
MatStencil 		rows,cols[5];

// specify matrix entries
vals[0]=4; 
vals[1]=-1; vals[2]=-1;
vals[3]=-1; vals[4]=-1;

// specify matrix rows and columns
row.i    = i  ;    row.j = j  ;
col[0].i = i  ; col[0].j = j  ;
col[1].i = i-1; col[1].j = j  ; 
col[2].i = i+1; col[2].j = j  ; 
col[3].i = i  ; col[3].j = j-1;
col[4].i = i  ; col[4].j = j+1;

// set the matrix values and entries
MatSetValuesStencil(A,1,&row,5,col,vals,INSERT_VALUES);

\end{lstlisting}

Because matrix \texttt{A} was associated with the DMDA object \texttt{da} at initialization, PETSc is able to map the node indices \texttt{(i,j)} to the corresponding matrix rows and columns. In addition, boundary conditions are more easily implemented using DMDA. For periodic boundary conditions, DMDA maps indices that are too large or small to the other edge of the domain. If, for example, \texttt{i-1=-1} in Listing~\ref{lst:dmda_entry}, then DMDA would map this entry to the matrix column corresponding to \texttt{i=nx-1}.

To solve the system $Ax=b$ in PETSc, we use the built in Krylov Subspace projection (KSP) methods with a preconditioner (PC). In this case, GMRES with a Jacobi preconditioner are used, but in the future we will implement user-defined preconditioners based on domain decomposition ideas. The PETSc \texttt{ksp} and \texttt{pc} objects are initialized independently of the DMDA object.

\begin{lstlisting}[label={lst:KSP},caption={Creating PETSc KSP and PC objects.}]
// set the KSP method
KSPCreate(PETSC_COMM_WORLD,&ksp);
KSPSetType(ksp,KSPGMRES);

// set the PC
KSPGetPC(ksp,&pc);
PCSetType(pc,PCJACOBI);
KSPSetPC(ksp,pc);
\end{lstlisting}

Finally, the system $Ax=b$ is solved with a call to \texttt{KSPSetOperators(ksp,A,A)} and \texttt{KSPSolve(ksp,b,x)}. However, we should note that special care is required for the all-Neumann boundary condition case for the Poisson equation~\ref{eq:pressuresolve}, which results in a singular system. This occurs, for example, when velocities are prescribed at all boundaries of the domain. To remedy this, we use PETSc to remove the null space of constant functions. 

\begin{lstlisting}[label={lst:nullspace},caption={Removing a nullspace in PETSc.}]
MatNullSpaceCreate(PETSC_COMM_WORLD,PETSC_TRUE,0,0,&nsp);
MatSetNullSpace(A,nsp);
MatNullSpaceRemove(nsp,b);
\end{lstlisting}

We evaluated the parallel performance of this implementation with a strong scaling test. The scaling and speedup relative to the single processor case are shown below for a test problem on a 256x256 and 512x512 grid. For the 256x256 grid, good scaling is seen up to 3 nodes, when parallel communication costs begin to degrade the performance. On the other hand, the 512x512 case has a better computation to communication ratio, and as such it scales nearly ideally out to 8 nodes.
\begin{figure}[!htb]
	\label{fig:scaling}
	\centering 
	\includegraphics[width=0.7\linewidth]{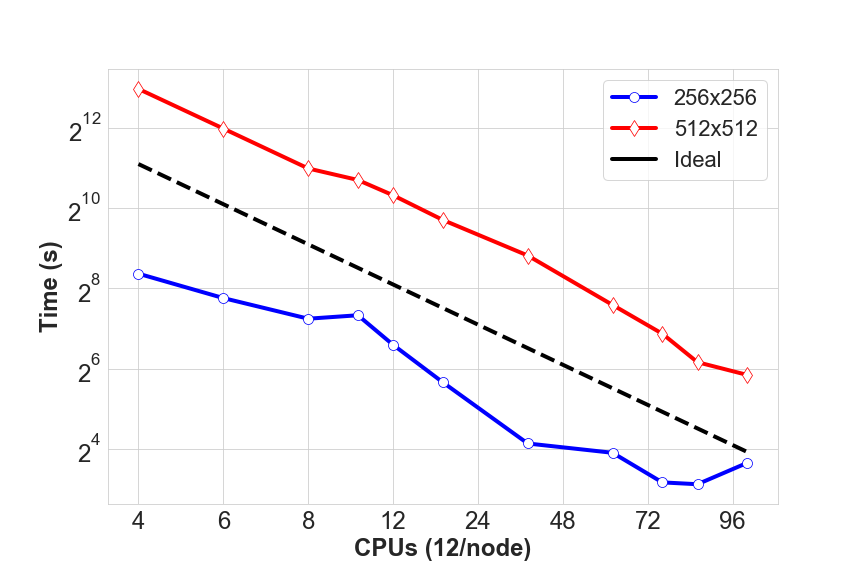}
	\caption{Strong scaling raw simulation times compared to an ideal scaling reference line.}
\end{figure}
\begin{figure}[!htb]
\label{fig:speedup}
	\centering 
	\includegraphics[width=0.7\linewidth]{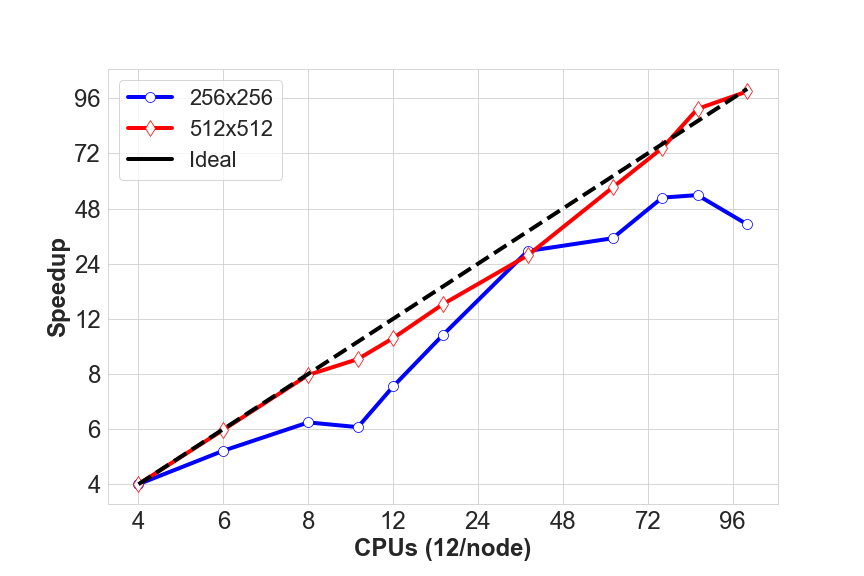}
	\caption{Strong scaling speedup compared to the ideal case.}
\end{figure}

\section{Numerical Results}\label{sec:NumericResults}
We simulated the evolution of the biofilm in a variety of settings.  For the result presented here we used periodic boundary conditions in the $x$ direction to be able to view downstream effects.  We used the dimensionless parameter values below for all cases unless otherwise noted.  
\begin{align}
&Re_s = 9.98\times 10^{-4} , \hspace{.25 in}Re_n = 2.33 \times 10^{-9}, \hspace{.25 in} \Lambda = 10^{-10}, \notag \\
&\Gamma_1 = 33.467, \hspace{.25 in}\Gamma_2 = 1.25\times 10^6,\hspace{.25 in} D_s = 2.3.  \\
& \mu = 0.14, \hspace{.25 in} K_c = 0.15, \hspace{.25 in} A = 100\notag
\end{align}
A table of all the dimensional parameters is included below.
\begin{table}[htbp]
\centering
\caption{List of parameters values used for the simulations.}
\label{tab:parameter4}
\begin{tabular}{||l|l|l|l||}
\hline
Symbol & Parameter & Value & Units\\
\hline
$T$ & Temperature & $303$ & K\\
$\gamma_1$ & Distortional energy & $5\times 10^7$ & $\text{kgm}^{-1}\text{s}^{-2}$\\
$\gamma_2$ & Separation energy & $1\times 10^{16}$ & $\text{kgm}^{-1}\text{s}^{-2}$\\
$\chi$ & Flory-Huggins parameter & 0.58 & \\
$\lambda_{CH}$ & Mobility parameter & $1\times 10^{-10}$ & $\text{kg}^{-1}m^3s$ \\
$N$ & Generalized polymerization parameter & $1 \times 10^3 $ & \\
$\mu$ & Max production rate & $1.4 \times 10^{-4}$ & $\text{kgm}^{-3}s^{-1}$ \\
$A$ & Max consumption rate & 0.1 & $\text{kgm}^{-3}s^{-1}$ \\
$D$ & Diffusion coefficient. & $2.3\times 10^{-9}$ & $\text{m}^2\text{s}^{-1}$ \\
$\eta_n$ & Dynamic viscosity of the network & $4.3\times 10^2$ & $\text{kgm}^{-1}s^{-1}$ \\
$\eta_s$ & Dynamic viscosity of the solvent & $1.002 \times 10^{-3} $ & $\text{kgm}^{-1}s^{-1}$ \\
$\rho_n$ & Network density & $1\times 10^3$ & $\text{kgm}^{-1}$ \\
$\rho_s$ & Solvent density & $1\times 10^3$ & $\text{kgm}^{-1}$ \\
$t_0$ & Characteristic time scale. & $1 \times 10^{3}$ & s \\
$h_0$ & Characteristic length scale & $1\times 10^{-3}$ & m \\
$c_0$ & Characteristic substrate concentration. & $1\times 10^{-3}$ & $\text{kgm}^{-3}$\\
$M$ & number of spatial intervals & 256 & \\
\hline
\end{tabular}
\end{table}
\subsection{2D Fluid Flow with Biofilm Protrusions.}
\begin{figure}[ht]
    \centering
   \includegraphics[width=0.7\textwidth]{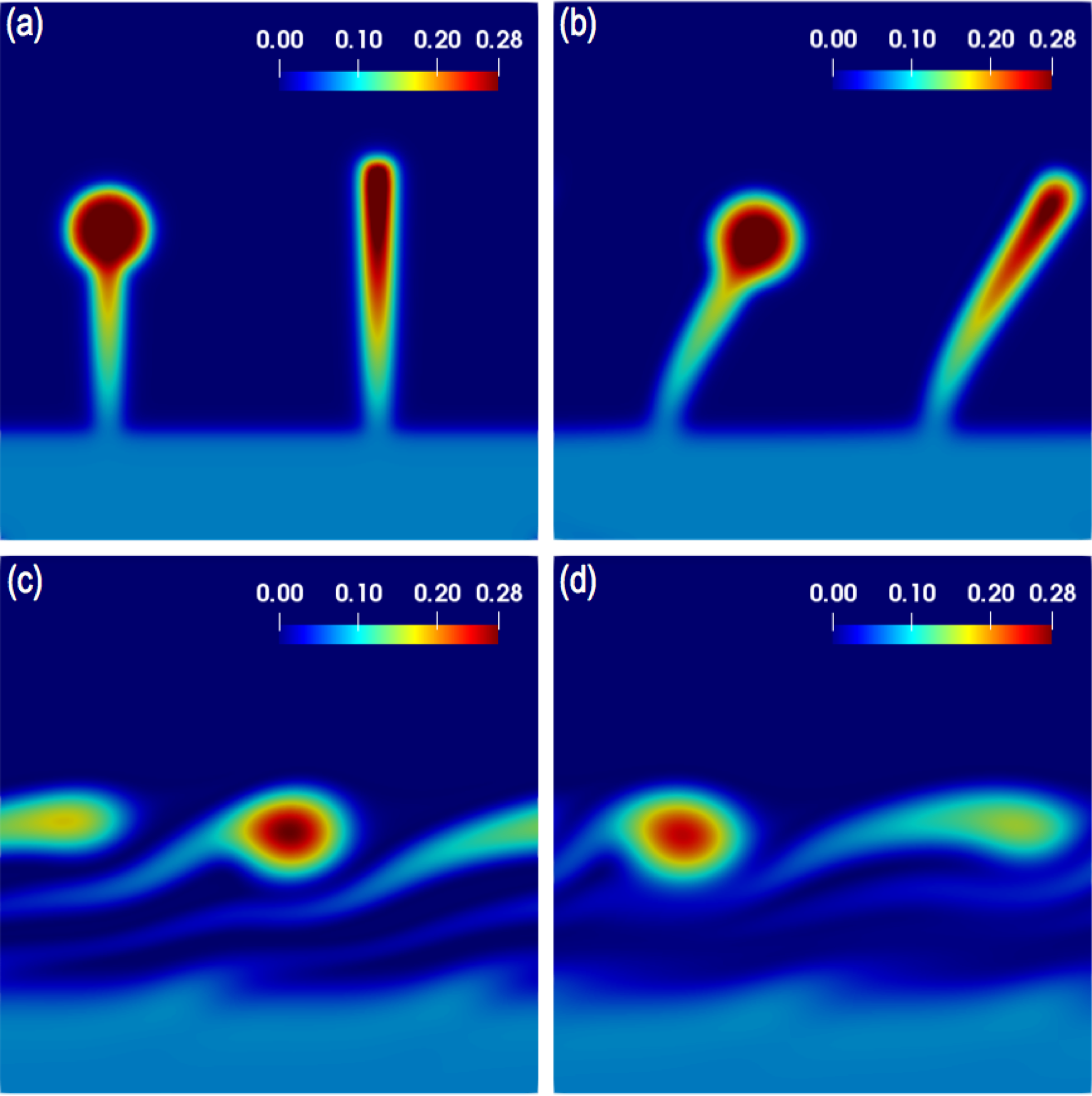}
  \caption[A cartoon of the biofilm in 2D]{The biofilm is now allowed to change in both the horizontal and vertical directions.  In part (a) we see the intial condition.  In part (b) we see the result after 20 timesteps or the equivilent of 5.5 hours.  Part (c) the initial condition has been stretched out and is starting to thin near the attachment point after 200 timesteps. Part (d), after timestep 500 equivilent to 5.8 days, the biofilm has broken away from the base layer of material and is now moving freely in the flow.}
  \label{fig:twofingerforms}
\end{figure}

We first look at a case to test the parallel code and compare it against previous results in the serial code.  In this case we start with a pair of mushroom shaped protrusions of biofilm away from the base.  There is a smaller concentration of biofilm in the neck region of the protrusion when compared to the top.  For this simulation we used periodic boundary conditions at the left and right boundaries and no flux of the polymer network at the top and bottom.  We also applied a shear velocity at the top of the domain \\
\begin{center}
\begin{equation}
u\big{|}_{y=1} = 0.1 \hspace{.25 in} v\big{|}_{y=1} = 0.
\end{equation}
\end{center}
As can be seen in figure \ref{fig:twofingerforms} the mushroom shaped regions are first deformed by the fluid flow.  Both regions undergo a stretching as well as a thining.  This is easier to see by comparing the large circular region at the top of the first protrusion across different timesteps.  Eventually the mushroom shaped regions break away from the base.  While this simulation is different from previous work it does exhibit the expected end behavior of shedding a region of biofilm that is connected through a thin neck to the bulk biofilm as described in \cite{Zhang_2D}.

\section{Conclusion and Future Work}
\label{sec:Conc}
We presented both the numerical methods used as well as the results of those methods.  The parallel implementation scales nearly ideally out to 8 nodes, and further scaling is expected for larger meshes. The simulations exhibit the expected behavior long term.  The protrusion of biofilm is connected by a thin neck to the bulk biofilm and eventually breaks away.  The pieces that have broken off then merge together into a single mass of biofilm and continue to move with the fluid flow present.  Work is currently underway on the combined Navier-Stokes-Cahn-Hilliard 3D model and a more generalized 2D model that will allow us to simulate more realistic geometries and flow fields.  We have also investigated implementing preconditioners based on domain decomposition methods  \cite{Kimn:2007, Kimn:2013, zheng2015parallel} in order to design an efficient parallel simulation procedure.

\bibliographystyle{siamplain}
\bibliography{MSDKbib}

\begin{thebibliography}{10}

\bibitem{alpkvista2007multidimensional}
{\sc E.~Alpkvist and I.~Klapper}, {\em A multidimensional multispecies
  continuum model for heterogeneous biofilm development}, Bulletin of
  mathematical biology, 69 (2007), pp.~765--789.

\bibitem{petsc-web-page}
{\sc S.~Balay, S.~Abhyankar, M.~F. Adams, J.~Brown, P.~Brune, K.~Buschelman,
  L.~Dalcin, A.~Dener, V.~Eijkhout, W.~D. Gropp, D.~Kaushik, M.~G. Knepley,
  D.~A. May, L.~C. McInnes, R.~T. Mills, T.~Munson, K.~Rupp, P.~Sanan, B.~F.
  Smith, S.~Zampini, H.~Zhang, and H.~Zhang}, {\em {PETS}c {W}eb page}.
\newblock \url{http://www.mcs.anl.gov/petsc}, 2018,
  \url{http://www.mcs.anl.gov/petsc}.

\bibitem{petsc-user-ref}
{\sc S.~Balay, S.~Abhyankar, M.~F. Adams, J.~Brown, P.~Brune, K.~Buschelman,
  L.~Dalcin, A.~Dener, V.~Eijkhout, W.~D. Gropp, D.~Kaushik, M.~G. Knepley,
  D.~A. May, L.~C. McInnes, R.~T. Mills, T.~Munson, K.~Rupp, P.~Sanan, B.~F.
  Smith, S.~Zampini, H.~Zhang, and H.~Zhang}, {\em {PETS}c users manual}, Tech.
  Report ANL-95/11 - Revision 3.10, Argonne National Laboratory, 2018,
  \url{http://www.mcs.anl.gov/petsc}.

\bibitem{chorin1967}
{\sc A.~J. Chorin}, {\em A numerical method for solving incompressible viscous
  flow problems}, Journal of computational physics, 2 (1967), pp.~12--26.

\bibitem{cogan2005channel}
{\sc N.~Cogan and J.~P. Keener}, {\em Channel formation in gels}, SIAM Journal
  on Applied Mathematics, 65 (2005), pp.~1839--1854.

\bibitem{costerton2003medical}
{\sc B.~Costerton}, {\em Medical biofilm microbiology: The role of microbial
  biofilms in disease}, Chronic Infections, and Medical Device Failure, CD-ROM,
  Montana State University,  (2003).

\bibitem{guermond_projectionmethod}
{\sc J.~Guermond, P.~Minev, and J.~Shen}, {\em An overview of projection
  methods for incompressible flows}, Computer methods in applied mechanics and
  engineering, 195 (2006), pp.~6011--6045.

\bibitem{Kimn:2007}
{\sc J.~Kimn and M.~Sarkis}, {\em Restricted overlapping balancing domain
  decomposition methods and restricted coarse problem for the {H}elmholtz
  equation}, Computer Methods in Applied Mechanics and Engineering, 196 (2007),
  pp.~1507--1514.

\bibitem{Kimn:2013}
{\sc J.~Kimn and M.~Sarkis}, {\em {Shifted Laplacian RAS solvers for Helmholtz
  equation}}, in Lecture Notes in Computational Science and Engineering,
  vol.~91, 2013, pp.~151--158.

\bibitem{klapper2002viscoelastic}
{\sc I.~Klapper, C.~Rupp, R.~Cargo, B.~Purvedorj, and P.~Stoodley}, {\em
  Viscoelastic fluid description of bacterial biofilm material properties},
  Biotechnology and Bioengineering, 80 (2002), pp.~289--296.

\bibitem{phillips2013engineered}
{\sc A.~J. Phillips, R.~Gerlach, E.~Lauchnor, A.~C. Mitchell, A.~B. Cunningham,
  and L.~Spangler}, {\em Engineered applications of ureolytic
  biomineralization: a review}, Biofouling, 29 (2013), pp.~715--733.

\bibitem{picioreanu1996modelling}
{\sc C.~Picioreanu}, {\em Modelling biofilms with cellular automata}, Final
  report to European Environmental Research Organisation,  (1996).

\bibitem{picioreanu2004particle}
{\sc C.~Picioreanu, J.-U. Kreft, and M.~C. van Loosdrecht}, {\em Particle-based
  multidimensional multispecies biofilm model}, Applied and environmental
  microbiology, 70 (2004), pp.~3024--3040.

\bibitem{picioreanu1999multidimensional}
{\sc C.~Picioreanu, M.~van Loosdrecht, and J.~Heijnen}, {\em Multidimensional
  modeling of biofilm structure}, Delft University of Technology, Faculty of
  Applied Sciences, 1999.

\bibitem{picioreanu1998mathematical}
{\sc C.~Picioreanu, M.~C. Van~Loosdrecht, J.~J. Heijnen, et~al.}, {\em
  Mathematical modeling of biofilm structure with a hybrid
  differential-discrete cellular automaton approach}, Biotechnology and
  bioengineering, 58 (1998), pp.~101--116.

\bibitem{picioreanu1998new}
{\sc C.~Picioreanu, M.~C. van Loosdrecht, J.~J. Heijnen, et~al.}, {\em A new
  combined differential-discrete cellular automaton approach for biofilm
  modeling: application for growth in gel beads}, Biotechnology and
  bioengineering, 57 (1998), pp.~718--731.

\bibitem{picioreanu2004advances}
{\sc C.~Picioreanu, J.~Xavier, and M.~C. van Loosdrecht}, {\em Advances in
  mathematical modeling of biofilm structure}, Biofilms, 1 (2004),
  pp.~337--349.

\bibitem{reiner1964deborah}
{\sc M.~Reiner}, {\em The deborah number}, Physics today, 17 (1964), p.~62.

\bibitem{spath1998sorption}
{\sc R.~Sp{\"a}th, H.-C. Flemming, and S.~Wuertz}, {\em Sorption properties of
  biofilms}, Water Science and Technology, 37 (1998), pp.~207--210.

\bibitem{tabak2005developments}
{\sc H.~H. Tabak, P.~Lens, E.~D. van Hullebusch, and W.~Dejonghe}, {\em
  Developments in bioremediation of soils and sediments polluted with metals
  and radionuclides--1. microbial processes and mechanisms affecting
  bioremediation of metal contamination and influencing metal toxicity and
  transport}, Reviews in Environmental Science and Bio/Technology, 4 (2005),
  pp.~115--156.

\bibitem{Zhang_2D}
{\sc T.~Zhang, N.~Cogan, and Q.~Wang}, {\em Phase field models for biofilms.
  ii. 2-d numerical simulations of biofilm-flow interaction}, Commun. Comput.
  Phys, 4 (2008), pp.~72--101.

\bibitem{Zhang_1D}
{\sc T.~Zhang, N.~G. Cogan, and Q.~Wang}, {\em Phase field models for biofilms.
  i. theory and one-dimensional simulations}, SIAM Journal on Applied
  Mathematics, 69 (2008), pp.~641--669.

\bibitem{zheng2015parallel}
{\sc X.~Zheng, C.~Yang, X.-C. Cai, and D.~Keyes}, {\em A parallel domain
  decomposition-based implicit method for the cahn--hilliard--cook phase-field
  equation in 3d}, Journal of Computational Physics, 285 (2015), pp.~55--70.

\end{thebibliography}
\end{document}